\documentclass[a4paper,12pt]{article}
\usepackage{amsmath, amsfonts, amssymb,xcolor}
\makeatletter
\newbox\bk@bxb
\newbox\bk@bxa
\newif\if@bkcont
\newcount\bk@lcnt

\def\breakboxskip{2pt}
\def\breakboxparindent{1.8em}

\def\breakbox{\vskip\breakboxskip\relax
\setbox\bk@bxb\vbox\bgroup
\advance\linewidth -2\fboxrule
\hsize\linewidth\@parboxrestore
\parindent\breakboxparindent\relax}

\def\bk@split{%
\@tempdimb\ht\bk@bxb 
\advance\@tempdimb\dp\bk@bxb
\setbox\bk@bxa\vsplit\bk@bxb to\z@ 
\setbox\bk@bxa\vbox{\unvbox\bk@bxa}
\setbox\@tempboxa\vbox{\copy\bk@bxa\copy\bk@bxb}
\advance\@tempdimb-\ht\@tempboxa
\advance\@tempdimb-\dp\@tempboxa}

\def\bk@addfsepht{%
\setbox\bk@bxa\vbox{\vskip\fboxsep\box\bk@bxa}}

\def\bk@addskipht{%
\setbox\bk@bxa\vbox{\vskip\@tempdimb\box\bk@bxa}}

\def\bk@addfsepdp{%
\@tempdima\dp\bk@bxa
\advance\@tempdima\fboxsep
\dp\bk@bxa\@tempdima}

\def\bk@addskipdp{%
\@tempdima\dp\bk@bxa
\advance\@tempdima\@tempdimb
\dp\bk@bxa\@tempdima}

\def\bk@line{%
\hbox to \linewidth{%
\hskip-2\fboxsep\vrule \@width\fboxrule\hskip.5\fboxsep\vrule \@width\fboxrule\hskip1.5\fboxsep
\box\bk@bxa\hfil
}}%

\def\endbreakbox{\egroup
\ifhmode\par\fi{\noindent\bk@lcnt\@ne
\@bkconttrue\baselineskip\z@\lineskiplimit\z@
\lineskip\z@\vfuzz\maxdimen
\bk@split\bk@addfsepht\bk@addskipdp
\ifvoid\bk@bxb 
\def\bk@fstln{\bk@addfsepdp
\hskip-\parindent\vbox{\llap{\raisebox{-2ex}{\rule{1.5\fboxsep}{\fboxrule}\hskip.5\fboxsep}}\bk@line\llap{\rule{1.5\fboxsep}{\fboxrule}\hskip.5\fboxsep}}}

\else 
\def\bk@fstln{\vbox{\llap{\raisebox{-2ex}{\rule{1.5\fboxsep}{\fboxrule}\hskip.5\fboxsep}}\bk@line}\hfil%
\advance\bk@lcnt\@ne
\loop
\bk@split\bk@addskipdp\leavevmode
\ifvoid\bk@bxb 
\@bkcontfalse\bk@addfsepdp
\vtop{\bk@line\noindent\hskip-2\fboxsep{\rule{1.5\fboxsep}{\fboxrule}}}%

\else 
\bk@line
\fi
\hfil\advance\bk@lcnt\@ne
\if@bkcont\repeat}%
\fi
\leavevmode\bk@fstln\par}\vskip\breakboxskip\relax}

\newcommand{\fracb}[2]{\frac{\raisebox{-.7ex}{$\scriptstyle #1$}}{\raisebox{.7ex}{$\scriptstyle #2$}}}

\def\smp{\smallskip\par}
\def\un{{\bf 1}}
\def\zero{\{0\}}
\def\pf{\noindent{\bf Proof~:}\ }
\def\findemo{~\leaders\hbox to 1em{\hss\  \hss}\hfill~\raisebox{.5ex}{\framebox[1ex]{}}\smp}

\def\mpn{\medskip\par\noindent}
\def\smpn{\smallskip\par\noindent}
\def\normal{\mathop{\trianglelefteq}}

\def\smpn{\smallskip\par\noindent}

\def\mpoint{\;\;.}
\def\mvirg{\;\;,}

\def\Res{{\rm Res}}

\def\Inf{{\rm Inf}}
\def\Def{{\rm Def}}

\def\Inf{{\rm Inf}}

\def\Z{\mathbb{Z}}
\def\N{\mathbb{N}}

\newcommand{\dirsum}[1]{\mathop{\oplus}_{#1}\limits}
\newcommand{\romain}[1]{\uppercase\expandafter{\romannumeral #1}}

\newcommand{\sur}[1]{\,\overline{\! #1}}
\newcommand{\sumb}[2]{\mathop{\sum}_{{\scriptstyle #1}\atop {\scriptstyle #2}}}

\newcommand{\sumc}[3]{\sum_{{\scriptstyle #1}\atop {{\scriptstyle #2}\atop {\scriptstyle #3}}}}

\newenvironment{enonce}[1]{\pagebreak[2]\refstepcounter{subsection}\refstepcounter{prop}\smpn{{\bf \thesection.\arabic{prop}.\ \ #1~:}}\begin{it} }{\end{it}\smp}
\newenvironment{enonce*}[1]{\pagebreak[2]\smpn{#1~:}\begin{it} }{\end{it}\smp}
\newcommand{\result}[1]{\begin{enonce}{#1}}
\def\fresult{\end{enonce}}
\newcommand{\npar}{\smallskip\par\noindent\pagebreak[2]\refstepcounter{subsection}\refstepcounter{prop}{\bf \thesection.\arabic{prop}.\ \ }}

\newenvironment{mth}[1]{\begin{breakbox}\begin{enonce}{#1}}{\end{enonce}\end{breakbox}}
\newenvironment{mth*}[1]{\begin{breakbox}\begin{enonce*}{#1}}{\end{enonce*}\end{breakbox}}
\newenvironment{rem}[1]{\refstepcounter{subsection}\refstepcounter{prop} \mpn{{\bf \thesection.\arabic{prop}.}\ \ \bf#1\ :}}{\smp}

\makeatletter
\renewenvironment{enumerate}{\ifnum \@enumdepth >3 \@toodeep\else
      \advance\@enumdepth \@ne
      \edef\@enumctr{enum\romannumeral\the\@enumdepth}\list
      {\csname label\@enumctr\endcsname}{\setlength{\topsep}{1ex}\setlength{\itemsep}{0pt}\usecounter
        {\@enumctr}\def\makelabel##1{\hss\llap{##1}}}\fi}{\endlist}
\renewenvironment{itemize}{\ifnum \@itemdepth >3 \@toodeep\else \advance\@itemdepth \@ne
\edef\@itemitem{labelitem\romannumeral\the\@itemdepth}%
\list{\csname\@itemitem\endcsname}{\setlength{\topsep}{1ex}\setlength{\itemsep}{0pt}\def\makelabel##1{\hss\llap{##1}}}\fi}
{\endlist}
\def\@sect#1#2#3#4#5#6[#7]#8{\ifnum #2>\c@secnumdepth
    \let\@svsec\@empty\else
    \refstepcounter{#1}\edef\@svsec{\csname the#1\endcsname .\hskip .5em}\fi
    \@tempskipa #5\relax
     \ifdim \@tempskipa>\z@
       \begingroup #6\relax
         \@hangfrom{\hskip #3\relax\@svsec}{\interlinepenalty \@M #8\par}%
       \endgroup
      \csname #1mark\endcsname{#7}\addcontentsline
        {toc}{#1}{\ifnum #2>\c@secnumdepth \else
                     \protect\numberline{\csname the#1\endcsname}\fi
                   #7}\else
       \def\@svsechd{#6\hskip #3\relax  
                  \@svsec #8\csname #1mark\endcsname
                     {#7}\addcontentsline
                          {toc}{#1}{\ifnum #2>\c@secnumdepth \else
                            \protect\numberline{\csname the#1\endcsname}\fi
                      #7}}\fi
    \@xsect{#5}}
\def\section{\@startsection {section}{1}{\z@}{-3.5ex plus-1ex minus
    -.2ex}{2.3ex plus.2ex}{\reset@font\Large\bf}}  

\makeatother

\def\mar[#1]{\ar@{-}[#1]|-{\object@{<}}}
\def\marb[#1]{\ar@{-}[#1]|{\object+{  }}}

\usepackage[all]{xy}

\begin{document}
\centerline{\Large\bf Fast decomposition of $p$-groups}\vspace{.2cm}\par
\centerline{\Large\bf in the Roquette category, for $p>2$}\vspace{.3cm}\par
\centerline{\bf Serge Bouc
}\vspace{.5cm}\par
{\footnotesize {\bf Abstract :} Let $p$ be a prime number. In \cite{roquette-category}, I introduced the {\em Roquette category} $\mathcal{R}_p$ of finite $p$-groups, which is an additive tensor category containing all finite $p$-groups among its objects. In $\mathcal{R}_p$, every finite $p$-group $P$ admits a canonical direct summand $\partial P$, called {\em the edge} of $P$. Moreover $P$ splits uniquely as a direct sum of edges of {\em Roquette $p$-groups}}.\par  
{\footnotesize In this note, I would like to describe a fast algorithm to obtain such a decomposition, when $p$ is odd.
}\vspace{1ex}\par
{\footnotesize {\bf AMS Subject classification :} 18B99, 19A22, 20C99, 20J15.\vspace{1ex}}\par
\section{Introduction}
Let $p$ be a prime number. The Roquette category $\mathcal{R}_p$ of finite $p$-groups, introduced in~\cite{roquette-category}, is an additive tensor category with the following properties~:
\begin{itemize}
\item Every finite $p$-group can be viewed as an object of $\mathcal{R}_p$. The tensor product of two finite $p$-groups $P$ and $Q$ in $\mathcal{R}_p$ is the direct product $P\times Q$.
\item In $\mathcal{R}_p$, any finite $p$-group has a direct summand $\partial P$, called {\em the edge} of~$P$, such that
$$P\cong\dirsum{N\normal P}\partial (P/N)\mpoint$$
Moreover, if the center of $P$ is not cyclic, then $\partial P=0$.
\item In $\mathcal{R}_p$, every finite $p$-group $P$ decomposes as a direct sum
$$P\cong\dirsum{R\in \mathcal{S}}\partial R\mvirg$$
where $\mathcal{S}$ is a finite sequence of {\em Roquette groups}, i.e. of $p$-groups of normal $p$-rank 1, and such a decomposition is essentially unique. Given the group~$P$, such a decomposition can be obtained explicitly from the knowledge of a {\em genetic basis} of $P$.
\item The tensor product $\partial P\times\partial Q$ of the edges of two Roquette $p$-groups $P$ and $Q$ is isomorphic to a direct sum of a certain number $\nu_{P,Q}$ of copies of the edge $\partial (P\diamond Q)$ of another Roquette group (where both $\nu_{P,Q}$ and $P\diamond Q$ are known explicitly.
\item The additive functors from $\mathcal{R}_p$ to the category of abelian groups are exactly the {\em rational $p$-biset functors} introduced in~\cite{bisetsections}. 
\end{itemize}
The latter is the main motivation for considering this category~: any structural result on $\mathcal{R}_p$ will provide for free some information on such rational functors for $p$-groups, e.g. the representation functors $R_K$, where $K$ is a field of characteristic 0 (see \cite{doublact}, \cite{fonctrq}, and L. Barker's article~\cite{rhetoric}), the functor of units of Burnside rings~(\cite{burnsideunits}), or the torsion part of the Dade group~(\cite{dadegroup}). \par
The decomposition of a finite $p$-group $P$ as a direct sum of edges of Roquette $p$-groups can be read from the knowledge of a genetic basis of $P$. The problem is that the computation of such a basis is rather slow, in general. For most purposes however, the full details encoded in a genetic basis are useless, and it would be enough to know the direct sum decomposition. \par
Hence it would be nice to have a fast algorithm taking any finite $p$-group $P$ as input, and giving its decomposition as direct sum of edges of Roquette groups in the category $\mathcal{R}_p$. This note is devoted to the description of such an algorithm, when $p>2$.
\section{Rational $p$-biset functors}
\npar Recall that the characteristic property of the edge $\partial P$ of a finite $p$-group in the Roquette category $\mathcal{R}_p$ is that for any rational $p$-biset functor $F$
$$\partial F(P)=\hat{F}(\partial P)\mvirg$$
where $\partial F(P)$ is the faithful part of $F(P)$, and $\hat{F}$ denotes the extension of $F$ to $\mathcal{R}_p$. Also recall the following criterion (\cite{rationnel}, Theorem 3.1):
\begin{mth}{Theorem} \label{rational}Let $p$ be a prime number, and $F$ be a $p$-biset functor. Then the following conditions are equivalent:
\begin{enumerate}
\item The functor $F$ is a rational $p$-biset functor.
\item For any finite $p$-group $P$, the following conditions hold:
\begin{itemize}
\item if the center of $P$ is non cyclic, then $\partial F(P)=\zero$.
\item if $E\normal P$ is a normal elementary abelian subgroup of rank 2, and if $Z\leq E$ is a central subgroup of order $p$ of $P$, then the map
$$\Res_{C_P(E)}^P\oplus \Def_{P/Z}^P:F(P)\to F\big(C_P(E)\big)\oplus F(P/Z)$$
is injective.
\end{itemize}
\end{enumerate}
\end{mth}
\npar Let $K$ be a commutative ring in which $p$ is invertible. When $P$ is a finite group, denote by $\mathsf{CF}_K(P)$ the $K$-module of central functions from $P$ to~$K$. The correspondence sending a finite $p$-group $P$ to $\mathsf{CF}_K(P)$ is a rational $p$-biset functor: 
\begin{mth}{Proposition} If $P$ and $Q$ are finite $p$-groups, if $U$ is a finite $(Q,P)$-biset, and if $f\in \mathsf{CF}_K(P)$, define a map $\mathsf{CF}_K(U):\mathsf{CF}_K(P)\to \mathsf{CF}_K(Q)$ by
$$\forall s\in Q,\;\; \mathsf{CF}_K(U)(f)(s)=\fracb{1}{|P|}\sumb{u\in U,\, x\in P}{su=ux}f(x)\mpoint$$
With this definition, the correspondence $P\mapsto \mathsf{CF}_K(P)$ becomes a rational $p$-biset functor, denoted by $\mathsf{CF}_K$.
\end{mth}
\pf A straightforward argument shows that $\mathsf{CF}_K(U)(f)$ is indeed a central function on $Q$, hence the map $\mathsf{CF}_K(U)$ is well defined. It is also clear that this map only depends on the isomorphism class of the biset $U$, and that for any two finite $(H,G)$-bisets $U$ and $U'$, we have
$$\mathsf{CF}_K(U\sqcup U')=\mathsf{CF}_K(U)+\mathsf{CF}_K(U')\mpoint$$
Moreover if $U$ is the identity biset at $P$, i.e. if $U=P$ with biset structure given by left and right multiplication, then for $f\in \mathsf{CF}_K(P)$ and $s\in P$
$$\mathsf{CF}_K(U)(f)(s)=\fracb{1}{|P|}\sumb{u\in U,\, x\in P}{su=ux}f(x)=\fracb{1}{|P|}\sum_{u\in P}f(s^u)=f(s)\mvirg$$
hence $\mathsf{CF}_K(U)$ is the identity map.\par
Now if $R$ is a third finite $p$-group, and $V$ is a finite $(R,Q)$-biset, then for any $t\in R$, setting $\lambda=\mathsf{CF}_K(V)\circ \mathsf{CF}_K(U)(f)(t)$, we have that
\begin{eqnarray*}
\lambda&=&\fracb{1}{|Q|}\sumb{v\in V,\,s\in Q}{tv=vs}\fracb{1}{|P|}\sumb{u\in U,\, x\in P}{su=ux}f(x)\\
&=&\fracb{1}{|Q||P|}\sumc{(v,u)\in V\times U}{s\in Q,\,x\in P}{tv=vs,\,su=ux} f(x)\\
&=&\fracb{1}{|Q||P|}\sumb{(v,u)\in V\times U,\,x\in P}{t(v,_{_Q}u)=(v,_{_Q}u)x}|\{s\in Q\mid tv=vs,\,su=ux\}|\,f(x)\\
\end{eqnarray*}
\begin{eqnarray*}
\lambda&=&\fracb{1}{|Q||P|}\sumb{(v,_{_Q}u)\in V\times_{_Q} U,\,x\in P}{t(v,_{_Q}u)=(v,_{_Q}u)x}|Q:Q_v\cap {_uP}||Q_v\cap {_uP}|\,f(x)\\
&=&\fracb{1}{|P|}\sumb{(v,_{_Q}u)\in V\times_{_Q} U,\,x\in P}{t(v,_{_Q}u)=(v,_{_Q}u)x}f(x)=\mathsf{CF}_K(V\times_QU)(f)(t)\mpoint
\end{eqnarray*}
Hence $\mathsf{CF}_K(V)\circ \mathsf{CF}_K(U)=\mathsf{CF}_K(V\times_QU)$, and $\mathsf{CF}_K$ is a $p$-biset functor.\par
To prove that this functor is rational, we use the criterion given by Theorem~\ref{rational}. Suppose first that the center $Z(P)$ of $P$ is non-cyclic. Let $E$ denote the subgroup of $Z(P)$ consisting of elements of order at most $p$. Then saying that $\partial \mathsf{CF}_K(P)=\zero$ amounts to saying that for any $f\in \mathsf{CF}_K(P)$, the sum
$$S=\sum_{Z\leq E}\mu(\un,Z)\Inf_{P/Z}^P\Def_{P/Z}^Pf$$
is equal to 0, where $\mu$ denotes the M\"obius function of the poset of subgroups of $P$ (or of $E$). Equivalently, for any $s\in P$
$$S(s)=\sum_{Z\leq E}\mu(\un,Z)\fracb{1}{|P|}\sumb{aZ\in P/Z,\,x\in P}{saZ=aZx}f(x)=0\mpoint$$
This also can be written as
\begin{eqnarray*}
S(s)&=&\sum_{Z\leq E}\mu(\un,Z)\fracb{1}{|P||Z|}\sumb{a\in P,\,x\in P}{saZ=aZx}f(x)\\
&=&\fracb{1}{|P|}\sum_{Z\leq E}\fracb{\mu(\un,Z)}{|Z|}\sum_{a\in P,\,z\in Z}f(s^a.z)\\
&=&\fracb{1}{|P|}\sum_{Z\leq E}\fracb{\mu(\un,Z)}{|Z|}\sum_{a\in P,\,z\in Z}f\big((sz)^a\big)\\
&=&\sum_{Z\leq E}\fracb{\mu(\un,Z)}{|Z|}\sum_{z\in Z}f(sz)\\
&=&\sum_{z\in E}\Big(\sum_{z\in Z\leq E}\fracb{\mu(\un,Z)}{|Z|}\Big)f(sz)\mpoint
\end{eqnarray*}
\begin{mth}{Lemma} Let $E$ be an elementary abelian $p$-group of rank at least 2. Then for any $z\in E$
$$\sum_{z\in Z\leq E}\fracb{\mu(\un,Z)}{|Z|}=0\mpoint$$
\end{mth}
\pf For $z\in E$, set $\sigma(z)=\sum_{z\in Z\leq E}\limits\fracb{\mu(\un,Z)}{|Z|}$. Assume first that $z\neq 1$, i.e. $|z|=p$. If $Z\ni z$ is elementary abelian of rank $r$, then $\mu(\un,Z)=(-1)^rp^{\binom{r}{2}}$, hence $\fracb{\mu(\un,Z)}{|Z|}=(-1)^rp^{\binom{r-1}{2}-1}=-\,\frac{1}{p}\,\mu(\un,Z/{<}z{>})$. Hence setting $\sur{Z}=Z/{<}z{>}$ and $\sur{E}=E/{<}z{>}$,
$$\sigma(z)=-\,\frac{1}{p}\sum_{\un\leq \sur{Z}\leq \sur{E}}\mu(\un,\sur{Z})=0\mvirg$$
since $|\sur{E}|>1$. Now
$$\sum_{z\in E}\sigma(z)=\sigma(1)+\sum_{e\in E-\{1\}}\sigma(z)=\sum_{z\in Z}\sum_{z\in Z\leq E}\fracb{\mu(\un,Z)}{|Z|}=\sum_{\un\leq Z\leq E}\mu(\un,Z)=0$$
hence $\sigma(1)=0$, completing the proof of the lemma.\findemo
It follows that $S(s)=0$, hence $S=0$, as was to be shown.\par
For the second condition of Theorem~\ref{rational}, suppose that $E$ is a normal elementary abelian subgroup of $P$ of rank 2, and that $Z$ is a central subgroup of $P$ of order $p$ contained in $E$. Let $f\in \mathsf{CF}_K(P)$ which restricts to 0 to $C_P(E)$, and such that
$$\forall sZ\in P/Z,\;\;(\Def_{P/Z}^Pf)(sZ)=\fracb{1}{p}\sum_{z\in Z}f(sz)=0\mpoint$$
Thus $f(s)=0$ if $s\in C_P(E)$. Assume that $s\notin C_P(E)$. Then for $e\in E$, the commutator $[s,e]$ lies in $Z$. Moreover the map $e\in E\mapsto [s,e]\in Z$ is surjective. it follows that for  any $z\in Z$, there exists $e\in E$ such that $s^e=sz$. Thus $f(sz)=f(s^e)=f(s)$. Hence $\Def_{P/Z}^Pf(s)=f(s)=0$. Hence $f=0$, as was to be shown.\findemo
\section{Action of $p$-adic units}
Let $\Z_p$ denote the ring of $p$-adic integers, i.e. the inverse limit of the rings $\Z/p^n\Z$, for $n\in\N-\zero$. The group of units $\Z_p^\times$ is the inverse limits of the unit groups $(\Z/p^n\Z)^\times$, and it acts on the functor $\mathsf{CF}_K$ in the following way: if $\zeta\in \Z_p^\times$ and $P$ is a finite $p$-group, choose an integer $r$ such that $p^r$ is a multiple of the exponent of $P$, and let $\zeta_{p^r}$ denote the component of $\zeta$ in $(\Z/p^r\Z)^\times$. For $f\in \mathsf{CF}_K(P)$, define $\widehat{\zeta}_P(f)\in \mathsf{CF}_K(P)$ by
$$\forall s\in P,\;\;\widehat{\zeta}_P(f)(s)=f(s^{\zeta_{p^r}})\mpoint$$
Then clearly $\widehat{\zeta}_P(f)$ only depends on $\zeta$, and this gives a well defined map
$$\widehat{\zeta}_P:\mathsf{CF}_K(P)\to \mathsf{CF}_K(P)\mpoint$$
One can check easily (see \cite{bisetfunctors} Proposition 7.2.4 for details) that if $Q$ is a finite $p$-group, and $U$ is a finite $(Q,P)$-biset, then the square
$$\xymatrix{
\mathsf{CF}_K(P)\ar[r]^-{\widehat{\zeta}_P}\ar[d]_-{\mathsf{CF}_K(U)}&\mathsf{CF}_K(P)\ar[d]^-{\mathsf{CF}_K(U)}\\
\mathsf{CF}_K(Q)\ar[r]^-{\widehat{\zeta}_Q}&\mathsf{CF}_K(Q)\\
}
$$
is commutative. In other words, we have an endomorphism $\widehat{\zeta}$ of the functor $\mathsf{CF}_K$. It is straightforward to check that for $\zeta, \zeta'\in \Z_p^\times$, we have $\widehat{\zeta \zeta'}=\widehat{\zeta}\circ\widehat{\zeta'}$, and that $\widehat{1}$ is the identity endomorphism of $\mathsf{CF}_K$. So this yields an action of the group $\Z_p^\times$ on $\mathsf{CF}_K$.\par
It follows in particular that when $n\in\N-\zero$, and $P$ is a finite $p$-group, if we set
$$F_n(P)=\{f\in \mathsf{CF}_K(P)\mid \forall s\in P,\,f(s^{1+p^n})=f(s)\}\mvirg$$
then the correspondence $P\mapsto F_n(P)$ is a subfunctor of $\mathsf{CF}_K$: indeed $F_n$ is the subfunctor of invariants by the element $1+p^n$ of $\Z_p^\times$.\par 
It follows that $F_n$ is a rational $p$-biset functor, for any $n\in \N-\zero$, hence it factors through the Roquette category $\mathcal{R}_p$. In particular, for any finite $p$-group $P$, if $P$ splits as a direct sum
$$P\cong \dirsum{R\in \mathcal{S}}\partial R$$
of edges of Roquette groups in $\mathcal{R}_p$, then there is an isomorphism
$$F_n(P)\cong \dirsum{R\in \mathcal{S}}\partial F_n(R)\mpoint$$
\begin{mth}{Notation} For a finite $p$-group $P$, and an integer $n\in \N-\zero$, let $l_n(P)$ denote the number of conjugacy classes of elements $s$ of $P$ such that $s^{1+p^n}$ is conjugate to $s$ in $P$. Also set $l_0(P)=1$.
\end{mth}
With this notation, for any finite $p$-group $P$, and any $n\in\N-\zero$, the $K$-module $F_n(P)$ is a free $K$-module of rank $l_n(P)$. In particular, if $P=C_{p^m}$ is cyclic of order $p^m$, then $F_n(P)$ has rank $l_n(P)=p^{\min(m,n)}$. Thus if $m>0$, then $\partial F_n(C_{p^m})$ has rank $p^{\min(m,n)}-p^{\min(m-1,n)}$, since $C_{p^m}\cong \partial C_{p^m}\oplus C_{p^{m-1}}$ in~$\mathcal{R}_p$.
\begin{mth}{Theorem} \label{main}Assume that a $p$-group $P$ splits as a direct sum
$$P\cong \un\oplus\dirsum{m=1}^\infty a_m\partial C_{p^m}$$
of edges of cyclic groups in the Roquette category $\mathcal{R}_p$, where $a_m\in\N$. Then 
$$\forall m\geq 1, \;\;a_m=\frac{l_m(P)-l_{m-1}(P)}{p^{m-1}(p-1)}\mpoint$$
\end{mth}
\pf For any $n \in\N-\zero$, we have
$$l_n(P)=1+\sum_{m=1}^\infty a_m(p^{\min(m,n)}-p^{\min(m-1,n)})=1+\sum_{m=1}^n a_m(p^{m}-p^{m-1})\mpoint$$
For $n\in\N-\zero$, this gives $l_n(P)-l_{n-1}(P)=a_n(p^n-p^{n-1})$.\findemo
\begin{mth}{Corollary} \label{decompose}Suppose $p>2$. If $P$ is a finite $p$-group, then
$$P\cong \un\oplus\dirsum{m=1}^\infty \fracb{l_m(P)-l_{m-1}(P)}{p^{m-1}(p-1)}\,\partial C_{p^m}$$
in the Roquette category $\mathcal{R}_p$.
\end{mth}
\pf Indeed for $p$ odd, all the Roquette $p$-groups are cyclic, hence the assumption of Theorem~\ref{main} holds for any $P$.\findemo
\setcounter{prop}{0}
\section*{Appendix}
\begin{rem}{A GAP function}
The following function for the GAP software (\cite{GAP4}) computes the decomposition of $p$-groups for $p>2$, using Corollary~\ref{decompose}:\par
\begin{footnotesize}
\begin{verbatim}
#
# Roquette decomposition of an odd order p-group g
# output is a list of pairs of the form [p^n,a_n]
# where a_n is the number of summands of g
# isomorphic to the edge of the cyclic group of order p^n
#
\end{verbatim}
\newpage
\begin{verbatim}
roquette_decomposition:=function(g)
local prem,cg,s,i,x,y,z,pn,u;
    if IsTrivial(g) then return [[1,1]];fi;
    prem:=PrimeDivisors(Size(g));
    if Length(prem)>1 then 
        Print("Error : the group must be a p-group\n"); 
        return fail;
    fi;
    prem:=prem[1];
    if prem=2 then 
        Print("Error : the order must be odd\n"); 
        return fail;
    fi;
    cg:=ConjugacyClasses(g);
    s:=[];
    for i in [2..Length(cg)] do
        x:=cg[i];
        y:=Representative(x);
        pn:=1;
        u:=y;
        repeat
            pn:=pn*prem;
            u:=u^prem;
            z:=y*u;
        until z in x;
        Add(s,pn);
    od;
    s:=Collected(s);
    s:=List(s,x->[x[1],x[2]*prem/(prem-1)/x[1]]);
    s:=Concatenation([[1,1]],s);
    return s;
end;
\end{verbatim}
\end{footnotesize}
\end{rem}
\begin{rem}{Example}
\begin{footnotesize}
\begin{verbatim}
gap> l:=AllGroups(81);;
gap> for g in l do
> Print(roquette_decomposition(g),"\n");
> od;
[ [ 1, 1 ], [ 3, 1 ], [ 9, 1 ], [ 27, 1 ], [ 81, 1 ] ]
[ [ 1, 1 ], [ 3, 4 ], [ 9, 12 ] ]
[ [ 1, 1 ], [ 3, 7 ], [ 9, 3 ] ]
[ [ 1, 1 ], [ 3, 7 ], [ 9, 3 ] ]
[ [ 1, 1 ], [ 3, 4 ], [ 9, 3 ], [ 27, 3 ] ]
[ [ 1, 1 ], [ 3, 4 ], [ 9, 4 ] ]
[ [ 1, 1 ], [ 3, 8 ] ]
[ [ 1, 1 ], [ 3, 5 ], [ 9, 1 ] ]
[ [ 1, 1 ], [ 3, 5 ], [ 9, 1 ] ]
[ [ 1, 1 ], [ 3, 5 ], [ 9, 1 ] ]
[ [ 1, 1 ], [ 3, 13 ], [ 9, 9 ] ]
[ [ 1, 1 ], [ 3, 16 ] ]
[ [ 1, 1 ], [ 3, 16 ] ]
[ [ 1, 1 ], [ 3, 13 ], [ 9, 1 ] ]
[ [ 1, 1 ], [ 3, 40 ] ]
\end{verbatim}
\end{footnotesize}
\end{rem}
For example, the group on line 6 of the previous list, isomorphic to the semidirect product $C_{27}\rtimes C_3$, is isomorphic to $\un\oplus4\partial C_3\oplus 4\partial C_9$ in $\mathcal{R}_3$.
\vspace{-2ex}

\centerline{\rule{5ex}{.1ex}}
\begin{flushleft}
Serge Bouc - CNRS-LAMFA, Universit\'e de Picardie, 33 rue St Leu, 80039, Amiens Cedex 01 - France. \\
{\tt email : serge.bouc@u-picardie.fr}\\
{\tt web~~ : http://www.lamfa.u-picardie.fr/bouc/}
\end{flushleft}
\end{document}